\documentclass[12pt]{amsart}

\usepackage{amssymb}
\usepackage{graphics}
\usepackage{graphicx}

\newtheorem{theorem}{Theorem}
\newtheorem{lemma}{Lemma}
\newtheorem{prop}{Proposition}
\newtheorem{defn}{Definition}
\newtheorem{cor}{Corollary}

\newcommand{\on}{\operatorname}
\newcommand{\PD}{\on{PD}}
\newcommand{\tb}{\on{tb}}
\newcommand{\id}{\on{id}}
\renewcommand{\d}{\partial}
\newcommand{\Spinc}{\on{Spin}^c}
\newcommand{\gr}{\ti{\on{gr}}}
\newcommand{\Fmix}{F^{\on{mix}}}
\newcommand{\F}{\mathcal{F}}
\newcommand{\D}{\mathcal{D}}

\newcommand{\Z}{\mathbb{Z}}
\newcommand{\R}{\mathbb{R}}
\newcommand{\Q}{\mathbb{Q}}
\newcommand\goth[1]{\mathfrak{#1}}
\newcommand{\s}{\goth{s}}
\renewcommand{\t}{\goth{t}}
\renewcommand{\k}{\goth{k}}
\newcommand{\di}{\on{d}}
\newcommand{\rk}{\on{rk}}
\newcommand{\ti}{\tilde}

\begin{document}

\author{Olga Plamenevskaya}  
\address{Department of Mathematics, Harvard University, Cambridge, MA 02138}
\email{olga@math.harvard.edu}
\title{Contact Structures with Distinct Heegaard Floer Invariants}

\begin{abstract}
We prove that the contact structures on $Y=\d X$ induced by non-homotopic
Stein structures on the 4-manifold $X$ have distinct Heegaard Floer 
invariants. 
\end{abstract}

\maketitle

\section{Introduction}
 
In \cite{LM}, Lisca and Mati\'c gave examples of 
non-isotopic contact structures which are homotopic as plane fields. 
Using Seiberg-Witten theory, they proved

\begin{theorem}\cite{LM} Let $W$ be a smooth compact 4-manifold with boundary, 
equipped with two Stein structures $J_1$, $J_2$ with associated $\Spinc$ 
structures $\s_1$, $\s_2$. If the induced contact structures on $\d W$ are
isotopic, then the $\Spinc$ structures $\s_1$ and $\s_2$ are isomorphic.  
\end{theorem}
  
In this paper we study Heegaard Floer contact invariants of such contact 
structures. These contact invariants were introduced by Ozsv\'ath and 
Szab\'o in \cite{ContOS}; to an oriented contact 3-manifold $(Y,\xi)$
with a co-oriented contact structure $\xi$ they 
associate an element $c(\xi)$ of the Heegaard Floer homology group
$\widehat{HF}(-Y)$. Conjecturally, Heegaard Floer homology is equivalent to 
Seiberg-Witten Floer homology, and the Heegaard Floer contact invariants 
are the same as the Seiberg-Witten invariants of contact structures 
constructed in \cite{KM}. In the Heegaard Floer context, we can make 
the theorem of Lisca and Mati\'c more precise:

\begin{theorem}\label{main}

Let $W$ be a smooth 4-manifold with boundary, 
equipped with two Stein structures $J_1$, $J_2$ with associated $\Spinc$ 
structures $\s_1$, $\s_2$, and let $\xi_1$, $\xi_2$ be the induced contact 
structures on $Y=\d W$. If the $\Spinc$ structures $\s_1$ and $\s_2$
are not isomorphic, then the contact invariants $c(\xi_1)$, $c(\xi_2)$ are 
distinct elements of $\widehat{HF}(-Y)$.   

\end{theorem} 
    
There is some additional structure on Heegaard Floer homology groups:
$\widehat{HF}(Y)$ decomposes as a direct sum 
$\bigoplus_{\s}\widehat{HF}(Y,\s)$ with summands
corresponding to $\Spinc$ structures on $Y$; if $c_1(\s)$ is torsion, 
the group $\widehat{HF}(Y, \s)$ is graded. 

If the contact structures $\xi_1$ and $\xi_2$ are homotopic as plane fields, 
they induce 
the same $\Spinc$ structure $\s$, and the contact invariants $c(\xi_1)$,
$c(\xi_2)$ both lie in $\widehat{HF}(-Y,\s)$. In the torsion case, 
they also have the same grading. However,  $c(\xi_1)$ and 
$c(\xi_2)$ can be nevertheless different, as follows from
Theorem \ref{main}.   

\medskip

{\bf Acknowledgements.} I am grateful to Peter Kronheimer for many helpful 
discussions, and to Paolo Lisca for pointing out a gap in the earlier 
version of this paper.

\section{Preliminaries on Heegaard Floer Homologies}

In this section we briefly recall some necessary facts
from the papers of Ozsv\'ath and Szab\'o \cite{3OS} - \cite{SympOS}. 

Given an oriented 3-manifold $Y$ equipped with a $\Spinc$ structure $\t$, 
there are homology groups $HF^+(Y,\t)$,
$HF^-(Y, \t)$, $\widehat{HF}(Y,\t)$. The last one is the simplest, but we will
mostly need the first two in this paper; the reader is referred 
to \cite{3OS}, \cite{3PrOS} for the definitions and properties. 
A cobordism between two 3-manifolds induces a map 
on homology. More precisely, if  $W$ is a cobordism from $Y_1$ to
$Y_2$, and $\s$ is a $\Spinc$ structure on $W$ with restrictions  
$\s|Y_1$, $\s|Y_2$ on $Y_1$, $Y_2$, then 
there are maps
$F_{W, \s}^{\circ}:HF^{\circ}(Y_1,\s|Y_1) \to HF^{\circ}(Y_2,\s|Y_2)$
($HF^\circ$ stands for one of the $HF^+$, $HF^-$, $\widehat{HF}$). 
These maps satisfy the composition law:

\begin{prop}\cite{4OS} Let $W_1$ be a cobordism from
  $Y_1$ to $Y_2$, $W_2$ a cobordism from $Y_2$ to $Y_3$, and 
$W=W_1\cup_{Y_2}W_2$ the composite cobordism. Let $\s_i\in\Spinc
  (W_i)$, $i=1,2$ be two $\Spinc$ structures with $\s_1|Y_1=\s_2|Y_2$.
Then for some choice of signs

\begin{equation}\label{complaw}
F_{W_2, \s_2}^\circ \circ F_{W_1, \s_1}^\circ = \sum_{\{\s\in \Spinc
 (W) :  \s_1|W_1=\s_1, \s_2|W_2=\s_2\}  } \pm F_{W, \s}^\circ.
\end{equation}

\end{prop}

For a cobordism $W$ from $Y_1$ to $Y_2$ with  $b_2^+(W)>1$, 
there is also a mixed invariant $\Fmix_{W,\s}: HF^-(Y_1, \s|Y_1)\to
HF^+(Y_2, \s|Y_2)$. It is defined by taking an ``admissible cut''
$N$, which separates $W$ into cobordisms $W_1$ from $Y_1$ to $N$ and $W_2$ 
from $N$ to $Y_2$ with $b_2^+(W_i)>0$, and composing 
$F_{W_1}^-:HF^-(Y_1)\to HF^-(N)$  and  $F_{W_2}^+:HF^+(N)\to HF^+(Y_2)$
in a certain way. 
We skip the details, as we can simply fix some admissible cut 
in our constructions. It follows from the composition law (\ref{complaw}) that 
\begin{equation}\label{compmix}
F_{W_2, \s_2}^+ \circ \Fmix_{W_1, \s_1} = \sum_{\{\s\in \Spinc
 (W) :  \s_1|W_1=\s_1, \s_2|W_2=\s_2\} } \pm \Fmix_{W, \s}.
\end{equation}

If $\s$ is a torsion $\Spinc$ structure, the homology groups
$HF^\circ(Y,\s)$
are graded; the grading takes values in $\Q$ and changes under 
cobordisms according to the following dimension formula.

\begin{prop} \cite{4OS} If $W$ is a cobordism from
  $Y_1$ to $Y_2$ endowed with a $\Spinc$ structure whose restriction 
to $Y_1$ and $Y_2$ is torsion, then 
\begin{equation}\label{dim} 
\gr(F_{W,\s}^\circ(\xi))-\gr(\xi)= \frac{c_1(\s)^2-2\chi(W)-3\sigma(W)}4 
\end{equation}
for any homogeneous element $\xi$.
\end{prop}
It follows that $\Fmix_{W,\s}$ affects the gradings in the same way.   

A closed 4-manifold  $X$ can be punctured in two points and regarded 
as a cobordism from $S^3$ to $S^3$; if $b_2^+(X)> 1$,
the mixed invariant of this cobordism
gives a closed manifold invariant $\Phi(X)$. If $X$ is symplectic,  
this invariant satisfies an important non-vanishing theorem
\cite{SympOS}. 
Below we state a version of this theorem for Lefschetz fibrations, 
rephrasing it in terms 
of mixed invariants for convenience. Abusing notation,
we denote by $X$  both the closed manifold and the corresponding 
cobordism from sphere to sphere. Recall that $HF^{\pm}(S^3)$ are given 
by 
$$
HF_k^-(S^3)=\Z \text{ in gradings } k<0, k \text{ even}; \ \ 
HF_k^+(S^3)=\Z,\  k\geq 0, k \text{ even}.    
$$         
  
\begin{theorem}\label{symp}\cite{SympOS} Let $\pi:X\to S^2$ be 
a relatively minimal 
Lefschetz fibration over the sphere with generic fiber $F$ of genus
$g>1$, and $b_2^+(X)>1$. Then for the canonical $\Spinc$ structure
$k$ the map $\Fmix_{X, k}$  sends the generator of $HF_{-2}^{-}(S^3)$ to 
the generator of 
$HF_0^{+}(S^3)$ (and vanishes on the rest of $HF^{-}(S^3)$).

For any other $\Spinc$ structure $\s\neq k$ with 
$\langle c_1(\s), [F]\rangle \leq 2-2g=\langle c_1(k), [F]\rangle $ 
the map $\Fmix_{X,\s}$ vanishes.

\end{theorem}

Finally, we need to recall the construction of the invariant
$c(\xi)$ for a contact manifold $(Y,\xi)$. 
We only consider co-oriented contact structures on $Y$. 
In \cite{ContOS}, $c(\xi)$ is defined as an element of $\widehat{HF}(-Y)$;
we will need to alter the definition slightly and look at  
$c(\xi)\in HF^+(-Y)$. The two elements obviously correspond to each other 
under the natural map $\widehat{HF}(-Y)\to HF^+(-Y)$. 
The definition uses the open book 
decomposion of $(Y,\xi)$, as well as  the following  
fact.

\begin{prop}\cite{SympOS} Let $Y_0$ be a fibration over the circle whose
fiber $F$ has genus  $g>1$. Let $\k$ be the canonical $\Spinc$ 
structure induced by the tangent planes 
to the fibers.  Then 
$$
HF^+(Y_0,\k)=\Z,
$$ 
and $HF^+(Y_0,\s)=0$ for any other $\s$ with 
$\langle c_1(\s), [F]\rangle =2-2g=\langle c_1(\k), [F]\rangle$.

\end{prop}

As proven by Giroux \cite{Gi}, contact manifolds can be described 
in terms of open books. An {\it open book decomposion} of $Y$
is a pair $(K, p)$  consisting of a (fibered) knot  $K\subset Y$
 and a fibration $p: Y\setminus K\to S^1$ whose fibers $p^{-1}(\phi)$ are 
interiors of compact embedded surfaces $F_\phi$ bounded by $K$. 
$K$ is then called the {\it binding} of the open book, and the fibers are
 the {\it pages}. An open book is compatible with 
a contact structure $\xi$ given by a contact form $\alpha$ on $Y$,  
if $\di\alpha$ is an area form on each page, and the binding is 
transverse to the contact planes and oriented as the boudary of 
$(F, \di\alpha)$. There is a one-to-one correspondence 
between isotopy classes of contact structures and the open books up to 
stabilization \cite{Gi}.

Given a compatible open book for $(Y, \xi)$, 
we can obtain a fibration  $Y_0$ by performing 0-surgery on the binding.
Let $V_0$ the corresponding cobordism from
$Y$ to $Y_0$, which can also be regarded as a cobordism from $-Y_0$
to $-Y$. The canonical $\Spinc$ structure $\k$ on $Y_0$ determines a
$\Spinc$ structure for $V_0$, so we can drop it from notation below.

\begin{defn} \cite[Proposition 3.1]{ContOS} The contact invariant 
is defined as 
$$
c(\xi)= F_{V_0}^+(c),  
$$ 
where $c$ stands for a generator of $HF^+(Y_0, \k)$, and 
$c(\xi)$ is defined up to sign.
\end{defn}

It is proven in \cite{ContOS} that $c(\xi)$ is independent of the 
choice of the open book. 

\section{Contact Invariants and Concave Fillings}

In this section we study contact invariants by using concave 
fillings of contact manifolds, and prove Theorem \ref{main}. 
More precisely, we prove  

\begin{theorem}\label{cobor} 
Let $W$ be a smooth compact 4-manifold with boundary $Y=\d X$. Let
$J_1$,  $J_2$ be two Stein structures on $W$ that induce $\Spinc$ structures 
$\s_1$, $\s_2$ on $W$ and contact structures  $\xi_1$, $\xi_2$ on $Y$. 
We puncture $W$ and regard it as a cobordism from
$-Y$ to $S^3$. Suppose that $\s_1|Y=\s_2|Y$, but the $\Spinc$ structures 
$\s_1$, $\s_2$ are not isomorphic. Then 
\begin{enumerate}
\item $F_{W, \s_i}^+(c(\xi_j))=0$ for $i\neq j$;

\item $F_{W,\s_i}^+(c(\xi_i))$ is a generator of $HF^+(S^3)$.
\end{enumerate}
\end{theorem}

Obviously, Theorem \ref{cobor} implies Theorem \ref{main} after 
we switch between the invariants in $\widehat{HF}(-Y)$ and their
images in  $HF^+(-Y)$: if $\s_1|Y\neq\s_2|Y$, then
the contact elements $c(\xi_1)$ and $c(\xi_2)$ lie in the components 
of $HF^+(-Y)$ corresponding different $\Spinc$-structures, 
and the statement  of Theorem \ref{main} is trivial.

\begin{cor} Suppose that a 3-manifold $Y$ bounds a compact smooth 
4-manifold $X$. If $X$ supports $n$ pairwise non-homotopic Stein 
structures, then 
$$
\rk(\widehat{HF}(Y))\geq n.
$$ 
\end{cor}
This follows from Theorem \ref{cobor} and the duality $\widehat{HF}^*(-Y)\cong
\widehat{HF}_*(Y)$, which gives an isomorphism $\widehat{HF}(-Y)\cong
\widehat{HF}(Y)$ in the non-torsion case \cite{3OS}.

For a Stein fillable contact manifold $(Y,\xi)$,
we want  to describe $c(\xi)$ as a mixed invariant of a certain 
concave filling of $(Y,\xi)$. 
We construct this concave filling,  following the work  
of Akbulut and Ozbagci \cite{AkO}. 

Suppose that $W$ is a (convex) Stein filling of $(Y, \xi)$. 
First we need to represent $W$  as a positive allowable Lefschetz 
fibration \cite{AO}, whose generic fiber is a surface with boundary. 
This induces  an open book decomposition
of $(Y, \xi)$ with monodromy consisting of non-separating positive Dehn 
twists. The original Stein structure on $W$ can be recovered from 
the Lefschetz fibration,  and the open book is compatible with $\xi$.
Note that the compatibility does not directly follow from the argument
in  \cite{AO}; for completeness we review this construction in Appendix,
strengthening it slightly and proving the compatibility statement.

Given a structure of positive allowable Lefschetz 
fibration on $W$ and the induced open  book decomposition of $Y$,
we perform $0$-surgery on the binding to get a cobordism $V_0$
from $Y$ to $Y_0$. Now  $W\cup  V_0$ is a Lefschetz fibration over the disk,
whose regular fiber $F$ is a closed surface obtained by capping off the page 
of the open book, and  $Y_0$ is a fibration over the circle.
The monodromy of $Y_0$ comes from the open book, and can be represented as 
a product $c_1 c_2 \dots c_k$ of positive Dehn twists. We need a concave 
filling for $Y_0$, so we want to construct a Lefschetz fibration 
with monodromy $c_k^{-1}\dots c_1^{-1}$. Recalling that the mapping class 
group of a closed surface is generated by non-separating positive 
Dehn twists \cite{AkO}, we rewrite  $c_k^{-1}\dots c_1^{-1}$
as a product of such twists. Putting in a node of 
the Lefschetz fibration for  each positive Dehn twist in the monodromy 
gives a  Lefschetz fibration $V_1$
with $\d V_1 =-Y_0$. 

\begin{figure}[ht]
\includegraphics[scale=1]{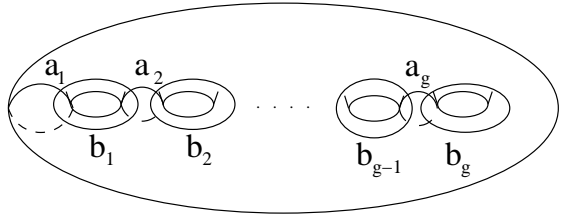} 
\caption{$(a_1 b_1 a_2 b_2\dots a_g b_g)^{4g+2}=1$}
\label{twists}
\end{figure}

To use the mixed invariants,
we must have $b_2^+(V_1)>1$. This can be achieved by a trick from 
\cite{AkO}: consider a Lefschetz fibration $G$ over the disk with the
nodes defined by Dehn twists of the word 
$(a_1 b_1 a_2 b_2\dots a_g b_g)^{4g+2}$; here  $a_i$, $b_i$ stand for 
positive Dehn twists around the curves shown in Figure \ref{twists}.
Since this word is equivalent to identity \cite{Bi}, we can glue 
a copy of $G$ into our fibration $V_1$ without affecting the monodromy of the
boundary $\d V_1=-Y_0$. On the other hand, a Lefschetz fibration 
has a symplectic structure \cite{GoSt}, so $b_2^+(G)>0$. Gluing in two copies 
of $G$, we may assume that $b_2^+(V_1)>1$.  

By construction, $X=W\cup V_0\cup V_1$ is a Lefschetz fibration over 
the sphere,
and $V=V_0\cup V_1$ is a concave symplectic filling of
$(Y, \xi)$.   
Let $k$ denote the canonical $\Spinc$ structure on $X$; we also write $k$
for its restrictions to $V$, $V_1$ etc.

Regard $V$ as a cobordism from $S^3$ to $-Y$, puncturing it at a point.       
The following fact is implicitly mentioned in \cite{ContOS}.
 
\begin{lemma}\label{concave} Suppose $c(\xi)$ is torsion. Let $\theta$ be the
  generator of $HF_{-2}^{-}(S^3)$. Then 
$$
c(\xi)=\pm\Fmix_{V, k}(\theta). 
$$ 
\end{lemma}     

\begin{proof} As before, let $c$ be the generator of $HF^+(-Y_0)$.
Observe that $\Fmix_{V_1}(\theta)=\pm c$. Indeed, by (\ref{compmix}) we have
 
$$
F_{V_0\cup W, k}^+\circ\Fmix_{V_1, k}=\sum_{\{\s\in \Spinc(X):\  
\  \s|V_0\cup W=k, 
\s |V_1=k \} }  \pm \Fmix_{X, s}=\pm \Fmix_{X, k}
$$
(there is just one term that survives in the sum, because $X$ is
a symplectic fibration, and the non-canonical $\Spinc$ structures
with $\langle c_1(\s), [F] \rangle=2-2g$ give 
zero maps by Theorem \ref{symp}). So  
$F_{V_0\cup W}^+\circ\Fmix_{V_1}(\theta)$ is the generator of $HF^+(S^3)$,
but this map factors through $HF^+(-Y_0,\k)=\Z$, so we must have   
$\Fmix_{V_1}(\theta)=\pm c$. 

Now it follows that   
\begin{equation}\label{sum}
c(\xi)=\pm F_{V_0}^+\circ\Fmix_{V_1, k}(\theta)=\sum_{\{\s\in\Spinc(V):\ 
 \s|V_0 =k, 
\s |V_1=k \}} \pm \Fmix_{V,\s}.
\end{equation}
Because $V_0$ consists of one $2$-handle attachment, 
the $\Spinc$-structures with given restrictions to $V_0$ and $V_1$ are 
of the form 
$k+n\PD[F]$, $n\in\Z$. The dimension formula (\ref{dim}) now 
implies that all non-zero terms in the sum (\ref{sum})
must have different absolute gradings, since 
$c_1(k+n\PD[F])^2=c_1(k)^2+2n(2-2g)$, and $g>1$. However, it is
clear from the definition that the contact invariant $c(\xi)$ is 
a homogeneous element in homology, so only one summand can be non-trivial.
This summand has to be  $\Fmix_{V, k}(\theta)$: again we can use the 
composition law (\ref{compmix}) and Theorem \ref{symp} to write 
\begin{equation}\label{nonvanish}
F_{W, k}^+\circ\Fmix_{V, k}=\sum_{\{\s\in \Spinc(X):\ \s|V_0\cup W=k, 
\s |V_1=k \} } \pm  \Fmix_{X, \s}=\pm \Fmix_{X, k},
\end{equation}
so $F_{W, k}^+\circ\Fmix_{V, k}(\theta)$ is the generator of $HF^+(S^3)$,
and it follows that $\Fmix_{V, k}(\theta)\neq 0$. Then
$\Fmix_{V, k}(\theta)=\pm c(\xi)$.
 
\end{proof}

\begin{proof}[Proof of Theorem \ref{cobor}]
We first deal with the case where $c_1(\s_i|Y)$ is torsion.  
Consider the Lefschetz fibration decomposition of 
the Stein manifold $(W, J_1)$,
and construct the concave filling $V$ as above  for the contact structure 
$\xi=\xi_1$. As before, the two pieces $V$ and $W$ fit together 
to form a Lefschetz
fibration $X$ over the sphere; by construction, $\s_1=k$ on $W$.
Looking at the proof of Lemma \ref{concave}, we can use (\ref{nonvanish}) 
to show that  $F_{W, \s_1}^+\circ\Fmix_{V, k}$ is the generator of 
$HF^+(S^3)$; by    Lemma \ref{concave} itself, 
$\Fmix_{V, k}(\theta)=\pm c(\xi_1)$, and Part (2) of the Theorem follows.

To prove Part (1), endow $W$ with the $\Spinc$ structure $\s_2$,
and glue it to the concave filling $V$ of  the contact structure 
$\xi=\xi_1$.  
Of course, we get the manifold $X$, which topologically remains the same, and 
the $\Spinc$ structures can be put together as $\s_1|Y=\s_2|Y$, but we  
no longer get the canonical $\Spinc$ structure associated to the symplectic 
structure.    
Again by Lemma \ref{concave} 
and the composition law,
$$
F_{W,\s_2}^+(c(\xi_1))=\pm F_{W,\s_2}^+\circ\Fmix_{V, k}(\theta)= 
\sum_{\{\s\in \Spinc(X):\  \s|W=\s_2, 
\s |V=k \} } \pm \Fmix_{X, \s}(\theta).
$$
For each of the $\Spinc$ structures in the sum we still have   
$\langle c_1(\s), [F] \rangle=2-2g$, but now none of them is canonical, 
since $\s_2$ is different from $\s_1$. By Theorem \ref{symp}, every term 
in the sum is zero. 

We have proved Theorem \ref{cobor} for the torsion case; 
it remains to treat the case when $c_1(\xi_i)$ is non-torsion. The
dimension formula is no longer valid and we can't use Lemma \ref{concave}, 
but we can look at the same construction and write
\begin{equation}\label{torsion}
F^+_{W,\s_i}(c(\xi_1))= \pm F_{W,\s_i}^{+}\circ F_{V_0, k}^+
\circ \Fmix_{V_1,k}(\theta)=\sum_\s \pm \Fmix_{X,\s},
\end{equation}
where the sum is now taken over all $\Spinc$ structures on $X$ 
which restrict to
$W$ as $\s_i$ and to $V_0$ and $V_1$ as $k$. As before, all these $\Spinc$
structures have $\langle c_1(\s), [F] \rangle=2-2g$.
If $\s_i=\s_1$,
one of the terms in (\ref{torsion}) corresponds to the canonical
$\Spinc$ structure on $X$, and the sum is equal to the generator of 
$HF^+(S^3)$; if $\s_i=\s_2$,  all the resulting $\Spinc$
structures on $X$  are different from
the canonical $\Spinc$ structure, so the sum is zero.  
\end{proof}

\section{An Example}
We now look at an example due to  Lisca and Mati\'c 
\cite{LM}.
\begin{figure}[ht]
\includegraphics[scale=1]{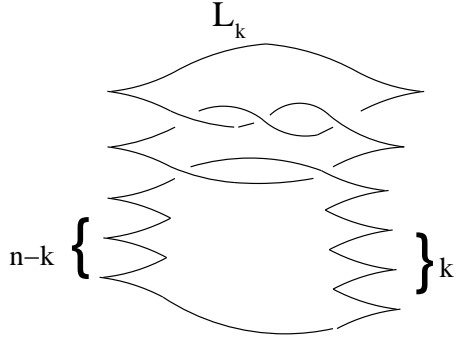} 
\caption{Legendrian link $L_k$}
\label{link}
\end{figure}

{\bf Example.} Let the contact manifold $(Y_n, \xi_k)$ be obtained as a 
Legendrian surgery on the Legendrian link $L_k$ shown on Fig. \ref{link} 
($k$ kinks on the right and $n-k$ kinks on the left  
give $r=2k-n$ for the rotation number of the unknot;
the rotation number of the trefoil is $0$). 
Varying $k$, we get  $n-1$ contact structures 
$\xi_1,\dots, \xi_{n-1}$ on $Y_n$. 

Topologically, the manifold $Y_n$ is the Brieskorn 
homology sphere $\Sigma(2, 3, 6n-1)$ with the orientation reversed; 
it is the boundary of the nucleus $N_n$. The Legendrian surgery cobordism
corresponding to $L_k$  endows $N_n$ with a Stein structure $J_k$. We have 
$c_1(J_k)=(2k-n)\PD[T]$, where $T$ is formed by a Seifert surface 
for the trefoil and the cocore of the handle attached to it \cite{LM},
so $J_k$ are pairwise non-homotopic. 
The contact  structures $\xi_i$ are all homotopic by Gompf's criterion
\cite{Go}, 
since  $Y$ is a homology sphere, and the Hopf invariant, defined as 
$h(\xi)=c_1(J)^2-2\chi(W)-3\sigma(W)$ for an almost-complex 4-manifold 
$(W, J)$ with boundary $(Y, \xi)$, is equal to $-6$ for all $\xi_i$. 
 
The manifold $Y$ can be obtained as $1/n$-surgery on the right-handed trefoil, 
and we can compute (cf. Section 8 of \cite{AbsOS})
$$
\widehat{HF}(-Y)=\widehat{HF}(-\Sigma(2, 3, 6n-1))=
\Z^n_{(+2)}\oplus\Z^{n-1}_{(+1)}, 
$$
where the subscripts indicate grading.

By \cite{ContOS}, 
the grading of the contact invariant is related to the Hopf invariant by 
$\gr(c(\xi))=-h(\xi)/4-1/2$, so for all $\xi_i$ the grading 
$\gr(c(\xi_i))=+1$. Theorem \ref{main} implies that 
the contact elements $c(\xi_i)$ 
are pairwise distinct; moreover, it follows from Theorem \ref{cobor} that
each $c(\xi_i)$ is a primitive element of $\widehat{HF}(-Y)$, and that 
$c(\xi_1), \dots c(\xi_{n-1})$ span
 $\Z^{n-1}_{(+1)}\subset \widehat{HF}(-Y)$.

{\it Remark.} The fact that the contact structures in this 
example have distinct contact invariants  was also discovered 
by Paolo Lisca and Andr\'as Stipsicz \cite{LiSt}, who have a 
different proof. 


\appendix{} 

\section{Lefschetz fibrations on Stein manifolds with boundary}

The decompositions of Stein manifolds as positive allowable 
Lefschetz fibrations
were constructed by Akbulut and Ozbagci in \cite{AO}. We give
an overview of their construction here, modifying it slightly 
to suit our purposes, and  taking particular care to prove that
 the open book induced by the Lefschetz 
fibration is compatible with the contact structure on the boundary 
of the Stein manifold.   

By \cite{El}, \cite{Go}, a Stein manifold $W$ with boundary $\d W=Y$ can be 
represented as $D^4 \cup(1\text{-handles})\cup(2\text{-handles})$; 
more precisely, 
$W$ is obtained by 
attaching $n$ $1$-handles to the ball $D^4$ and extending the 
Stein structure on $D^4$
to the handles to get the (unique) Stein structure on $\#_n S^1\times S^2$;
the $2$-handles are attached to components $L_i$ of a Legendrian link $L$ in 
$\#_n S^1\times S^2$, with the framings given by $\tb(L_i)-1$ ($\tb(L_i)$ 
denotes the Thurston-Bennequin number of $L_i$). In other words, $Y$
can be obtained by a Legendrian surgery on $L\subset \#_n S^1\times S^2$,
so that the corresponding surgery cobordism is $W$. 

To construct a Lefschetz fibration, we start with the case 
where $W$ has no $1$-handles, 
so $Y$ is obtained as a Legendrian surgery on a Legendrian link $L$ in $S^3$.

Let $\xi_0$ denote the standard contact structure on $S^3$.
The key ingredient of the construction is the following fact.
\begin{prop}\label{linkinbook} Given a Legendrian link $L\subset S^3$, 
there exist an open book 
decomposition of $S^3$, such that: 
\begin{enumerate}
\item
 the induced contact structure $\xi$ is isotopic 
to $\xi_0$;
\item 
   the link $L$ is contained in one of the pages, and does not separate it;  
\item 
$L$ is Legendrian with respect to $\xi$;  
\item 
there exist an isotopy which fixes $L$ and takes $\xi$ to $\xi_0$, so
the Legendrian type of 
the link is the same with respect to $\xi$ and $\xi_0$;
\item the framing of $L$ given by the page of the open book is the 
same as the contact framing.
\end{enumerate}    
\end{prop} 
Note that (5) trivially follows from (1)-(4).

 In \cite{AO} the statements (1), (2), and (5) are proved 
by putting $L$ into a ``square bridge 
position'' and constructing an open book for $S^3$ which contains this link.
The binding of this open book is a torus knot, so the monodromy 
produces the standard contact structure on $S^3$. Note, however, that
when $L$ is moved to the square bridge position, its Legendrian type is lost,
and the contact structure forgotten; 
Parts (3) and (4) of Proposition \ref{linkinbook} are unclear from \cite{AO}.

\begin{proof}
We start by constructing 
one page of the open book.  
Puncturing the sphere at one point,
we may consider links in $(\R^3,\xi_0)$; we assume that the contact 
structure  
$\xi_0$ on $\R^3$  is given by the contact form 
$\alpha_0= \di z+x \di y$. The next lemma is very similar to 
Theorem 2 from \cite{AO}, but keeps the link Legendrian and remembers 
the contact structure.   

\begin{lemma}\label{page} Given a Legendrian link $L$ in $(\R^3, \xi_0)$,
there exists a surface $F$ containing $L$, such that  $\di \alpha_0$
 is an area form for $F$, $\d F=K$ 
is a torus knot transverse to $\xi_0$, and $L$ does not separate $F$.  
\end{lemma}        
\begin{proof} After an appropriate Legendrian isotopy, we assume 
that the front projection of $L$ consists of segments which are straight lines 
(except in the neighborhoods of junctions), and  
all the negatively sloped segments have slope $-1$, while the positive slopes 
are all equal to $+1$ (see Figure \ref{sqleg}). Let $l_i=\{x=1, z=y+b_i\}$,
$i=1,\dots, p$ and  $m_j=\{x=-1, z=-y+d_j\}$,
$j=1,\dots, q$ be the lines in $\R^3$ containing  these segments; adding 
some extra lines if necessary, we can take $p$ and $q$ relatively 
prime. Denote by $t_{i,j}$ the intersection point of the lines $\{ z=y+b_i\}$
and  $\{ z=-y+d_j\}$ on the $yz$-plane.
\begin{figure}[ht]
\includegraphics[scale=0.5]{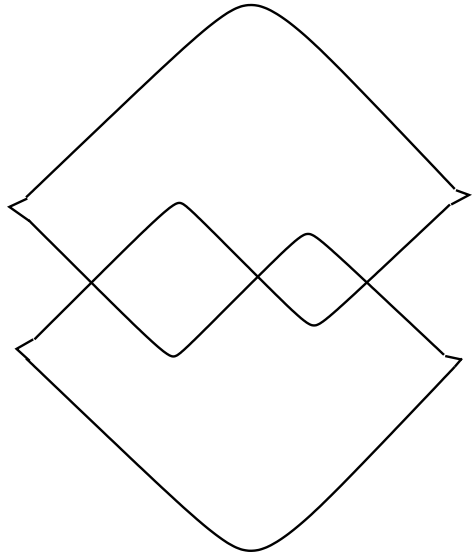} 
\caption{Legendrian link $L$}
\label{sqleg}
\end{figure}

We start the contruction of $F$ by looking at the narrow strips 
\begin{eqnarray*}
S^+_i&=&\{ 1-\epsilon \leq x\leq 1+\epsilon, z=y+b_i \}, i=1,\dots, p, \\   
S^-_j&=&\{ 1-\epsilon \leq x\leq 1+\epsilon, z=-y+d_j \}, j=1,\dots, q, 
\end{eqnarray*}
surrounding the straight segments of knots. Taking the strips long 
enough, so that each of the points $(\pm 1, t_{i,j})$ 
is contained in one of the strips,  we obtain a grid similar to the one 
shown on Figure \ref{grid}. If $\epsilon>0$ is small enough, 
$\di \alpha$ gives an area form on each strip. 
\begin{figure}[ht]
\includegraphics[scale=0.8]{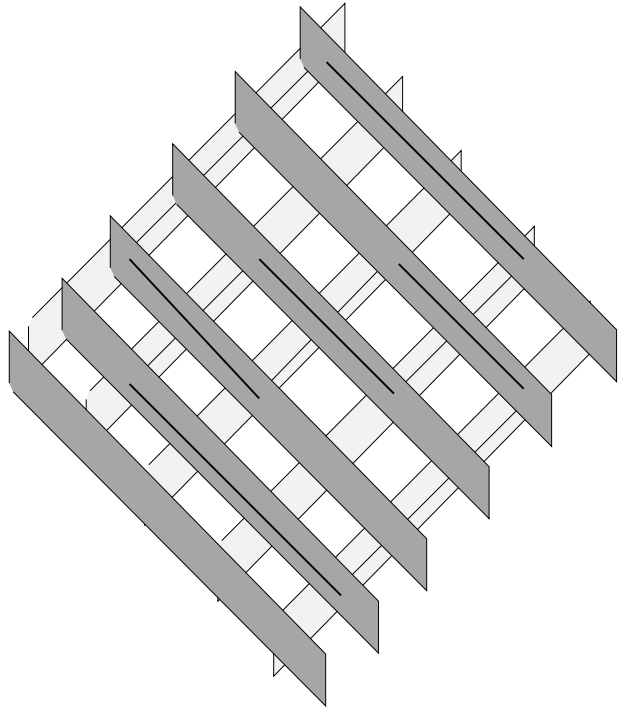} 
\caption{Strips of surface around L.}
\label{grid}
\end{figure}
\begin{figure}[ht]
\includegraphics[scale=0.7]{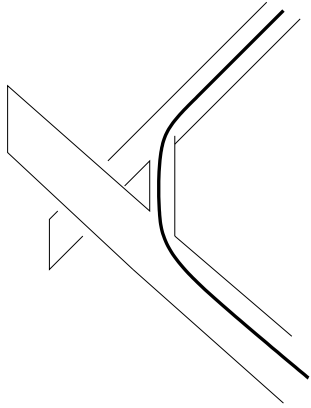} 
\caption{A connecting band.}
\label{band}
\end{figure}
We connect the points $(1, t_{i,j})$ and  $(-1, t_{i,j})$ for all $i,j$
by a segment of a straight line, and construct a band around this segment.
The band connects the strips $S^+_i$ and $S^-_j$; it  twists by 
$90^\circ$ along the way, following the contact planes 
(see Figure \ref{band}). By construction, we get a surface $F$ whose 
boundary is a $(p,q)$-torus knot; it is also clear that $F$ is close enough 
to the contact planes, which means that $\di \alpha_0$ induces 
an area form on $F$.

The line segments of the Legendrian link $L$ lie on $F$, but $L$ might 
not be contained in $F$ around the junctions. However, we can perturb $F$ 
slightly, and  move $L$ by a Legendrian isotopy to put it on $F$;
obviously, $L$ is non-separating.

It remains to ensure that the boundary of $F$ is transverse to the contact
planes. This is easy to achieve by moving the torus knot $K$ on $F$ to make 
it transverse to the characterictic foliation on $F$ (we may assume 
that the singular points of this foliation are isolated, so they do not 
present a problem).             
\end{proof}  

Now we can construct the required open book, starting with the page from 
Lemma \ref{page}.
We immediately get an interval worth of pages, 
perturbing this page slightly 
and making sure that $\di \alpha_0$ 
is an area form on each newly constructed page. 
The pages will span a handlebody 
$H_1$  (a thickening of the original page).   

Since  the torus knot $K$ is fibered, we can fiber the complementary 
handlebody $H_2=S^3\setminus H_1$ by the pages with binding $K$, 
thus completing the picture to a 
fibration $\pi:S^3\setminus K \to S^1$. 
Unfortunately,  the resulting open book
does not have to be compatible with the 
contact form $\alpha_0$: 
we have no guarantee that $\alpha_0$ induces 
an area form on pages in $H_2$. 
 
However, we can find a  contact form which is compatible 
with the open book $(K,\pi)$
and restricts to $\alpha_0$ on $H_1$ 
by using Thurston-Winkelnkemper 
construction \cite{ThWi}.  
Denote by $\phi:F\to F$ the monodromy of the open book, assuming that
$\phi=\id$ in the neighborhood of the binding. Let $T_\phi$ be the mapping 
torus 
$$
T_\phi=(F\times [0,1]) \slash \thicksim, 
\text{ where } (x,1)\thicksim(\phi(x),0). 
$$    
Cut out a small tubular neighborhood $K\times D^2$ 
of $K$ and shrink $F$ accordingly to represent the sphere as 
$$
S^3= T_\phi\cup_\d K\times D^2. 
$$
We may assume that the handlebody $H_1$ consists of the pages $F\times\{t\}$
with $t\in[0, 1/2]\subset S^1$, and that $\di\alpha_0$ gives an area form 
for all pages in the bigger handlebody $F\times [-\delta, 1/2+\delta]$ 
for some small $\delta>0$.
  
Let $\alpha^t$
be the restriction of the form $\alpha_0=\di z +x \di y$ to the page
$F\times\{t\}$. Set $\beta^t= 2(1-t)\alpha^{1/2}+ (2t-1)\phi^*\alpha^0$,
 $t\in [1/2,1]$; 
if 
$\kappa>0$ is large enough, $\beta^t+ \kappa \di t$ is a contact form on
$F\times [1/2,1]$. 

Let $\nu(t)$ be a positive increasing smooth function on $[1/2, 1/2+\delta]$,
such that $\nu(1/2)=0$, and $\nu(1/2+\delta)=\kappa$. The form 
$\alpha^t+\nu(t)\di t$ is contact on $F\times[1/2, 1/2+\delta]$, and 
``connects''   $\alpha_0$ and $\beta^t+ \kappa \di t$. Construct
$\alpha^t+\mu(t)\di t$ on $F\times[-\delta, 0]$ by analogy; now the forms
\begin{eqnarray*}
\alpha_0 &\text{ on }& F\times [0,1/2], \\
 \alpha^t+\nu(t)\di t &\text{ on }& F\times[1/2, 1/2+\delta], \\
\beta^t+ \kappa \di t &\text{ on }& F\times [1/2+\delta, 1-\delta], \\
\alpha^t+\mu(t)\di t &\text{ on }& F\times[-\delta, 0]
\end{eqnarray*}
fit together to produce a contact form on the mapping torus. This form 
extends over the binding, since the pages are 
transverse to the contact planes along the boundary.

As the binding of the open book $(K, \pi)$ is a torus knot, the monodromy 
is a product of non-separating positive Dehn twists \cite{AO}, 
so the corresponding 
contact structure $\xi$ is Stein fillable. It follows that 
$\xi$ is  isotopic to $\xi_0$. Moreover, the restrictions of $\xi$
and $\xi_0$ to the handlebody $H_1$ coincide, so the link $L$ obviously
remains Legendrian for $\xi$, and Part (3) of the 
Proposition is established.
To prove Part (4),  we will show that on the handlebody $H_2$ 
the contact structures $\xi_0$ and $\xi$ are isotopic relative to the 
boundary (as they coincide on $H_1$, the restrictions of $\xi$ and $\xi_0$
to $\d H_2=\d H_1$ are the same). We will be using convex 
surfaces and dividing curves (see \cite{Gir}, \cite{Ho}, \cite{Ka})
in our proof. 
Note that we may perturb the surface 
$\d H_1=\d H_2$ slightly, and assume that it is convex.

\begin{lemma}\label{handlebody} Assume that the 
handlebody $H\subset S^3$ is   
a thickening of the Seifert surface of a torus knot $K$. Consider tight 
contact structures on $H$ with convex boundary $\d H$, for which 
the dividing set is $\Gamma=K$. Suppose that two tight contact structures
$\xi$, $\zeta$ on $H$ induce the same  characteristic foliation
 $\F$ on $\d H$, and that $\F$ is adapted to $\Gamma$. Then 
$\xi$ and $\zeta$ are isotopic
 relative to $\d H$.     
\end{lemma}

\begin{proof} A Seifert surface of the (p,q)-torus knot $K$ can be obtained
by plumbing together $pq$ positive Hopf bands. The handlebody $H$ 
then decomposes as a boundary connected sum of thickened Hopf bands,
which can be thought of as solid tori with dividing set $\Gamma$ given 
by two parallel curves with slope $-1$. For any handlebody 
which is obtained by thickening of a plumbed sum of $n$ positive Hopf 
bands, we prove the statement of Lemma \ref{handlebody} by induction  
on $n$. The base of induction follows from 
Honda's classification of tight contact structures on the solid torus.

\begin{lemma}\cite{Ho} Let $\Gamma\subset T^2=\d S^1\times D^2$
consist of two parallel curves with slope $-1$. Then two tight  
contact structures on $S^1\times D^2$ with convex boundary $T^2$
are isotopic $\text{rel } T^2$ if they induce the same characterictic 
foliation 
adapted to $\Gamma$.   
\end{lemma}
   
For the induction step, we want to cut one of the solid tori 
off the handlebody $H$. Choose an appropriate disk $\D\subset H$
with convex boundary $\d\D\subset \d H$, so that $H=(S^1\times D^2)\cup_{\D} 
\ti{H}$, where $\ti{H}$ is a handlebody of smaller genus. 
To check that the dividing set is given by 
two curves of slope $-1$ on the boundary $S^1\times D^2$ and satisfies
 our assumption on the boundary of $\ti{H}$,
 we examine the dividing curves on the cutting disk $\D$.  
Observe that $\d \D$  meets $\Gamma$ 
in four points; 
we claim that after we cut along $\D$, round the corners of the resulting 
surfaces,  and regard 
$\D$ as part of $\d(S^1\times D^2)$ or $\d H$,  
the dividing set inside $\D$ consists 
of two lines joining these points pairwise (for each surface). 
Indeed, otherwise $\Gamma$ would 
have a component bounding a disk inside $\D$, but this is a contradiction 
with Giroux's criterion:
\begin{prop}(Giroux)  If $\Sigma\neq S^2$ is a convex surface (closed
or compact with Legendrian boundary) in a contact manifold $(M, \xi)$,
then $\Sigma$ has a tight neighborhood if and only if  the dividing set
$\Gamma_\Sigma$ has no homotopically trivial curves. 
\end{prop} 

\begin{figure}[ht]
\includegraphics[scale=0.5]{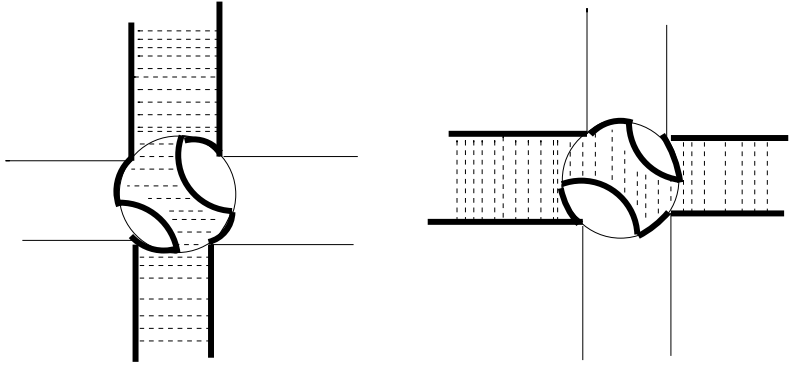} 
\caption{Correct dividing set.}
\label{goodcut}
\end{figure}

\begin{figure}[ht]
\includegraphics[scale=0.5]{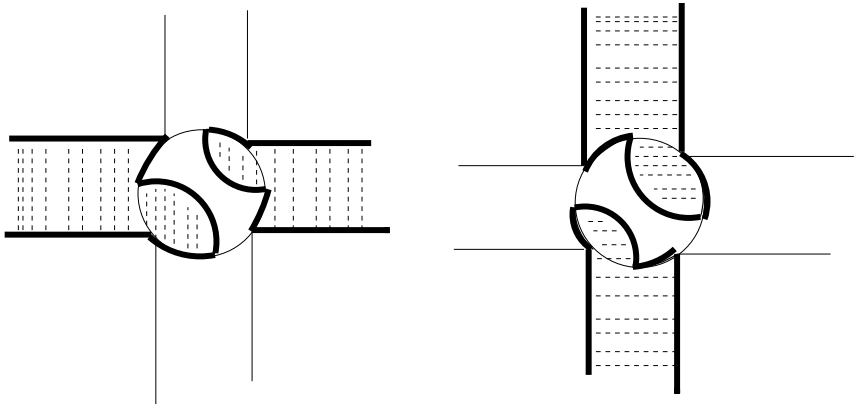} 
\caption{Wrong dividing set.}
\label{badcut}
\end{figure}

There are two possible ways for the two lines to join four points;
these are shown on Figures \ref{goodcut} and \ref{badcut} 
(for the explanation of what happens
to the dividing set when two convex surfaces meet and the corners are 
rounded, we refer the reader to  Section 3.3.2 in \cite{Ho}).
It remains to observe that  Figure \ref{badcut}  produces a
 homotopically trivial curve on the 
boundary of the solid torus and is ruled out by Giroux's criterion,
while on Figure \ref{goodcut} the dividing set 
on $\D$ connects the bands as required, decomposing the ``core surface'' 
of $H$ into a plumbed sum of a Hopf band and a ``core surface'' for $\ti{H}$.
The uniqueness of the tight contact structures on $\ti{H}$ and  
$S^1\times D^2$  (with given boundary conditions) now implies the uniqueness
of the tight contact structure on $H$, and the induction step follows.
(The dotted lines on Figures \ref{goodcut} and \ref{badcut}  are
used to highlight the bands and do not encode any foliation).
\end{proof}

The proof of Proposition \ref{linkinbook} is complete.
\end{proof}

Returning to the Lefschetz fibration construction, 
we can now obtain a required   
decomposition of a Stein manifold $W$ without $1$-handles: we represent 
$W$ as a Legendrian surgery cobordism for a Legendrian link $L$,
use Proposition \ref{linkinbook} to find an appropriate open book,
and add to the fibration a Lefschetz handle corresponding to 
the positive Dehn twist along a component of $L$ for each Legendrian 
2-handle of $W$ (see \cite{AO}); note that Lefschetz fibration given 
by the torus knot 
(in the absence of L) produces the (unique) Stein structure on $D^4$.

For the case where $1$-handles are present, we combine the argument 
from \cite{AO} with Proposition  \ref{linkinbook}. The Stein manifold
is represented as a Legendrian surgery on a link in $\#_n S^1\times S^2$,
which in turn corresponds to a diagram consisting of a Legendrian link 
in $S^3$ and $n$ dotted circles for the $1$-handles. We first use
Proposition  \ref{linkinbook} to find a ``nice'' open book for $S^3$,
and then for each dotted circle we scoop a disk out of each page,
so that the open book represents $\#_n S^1\times S^2$ now, 
and the Legendrian link is contained in a page. As before, we add 
Lefschetz handles compatibly with Legendrian handles. The pages 
of resulting open book will have multiple boundary components;
we need an open book with a connected binding to use the 
Ozsv\'ath-Szab\'o definition of contact invariants, so we make 
the boundary of the page connected by plumbing in some positive Hopf bands
(for the Stein fillings, this corresponds to taking the boundary 
connected sum with  the Stein ball $D^4$).

\end{document}